\documentclass[zamm,a4paper,fleqn]{w-art}
\usepackage{times}
\usepackage{w-thm}
\usepackage{cite}
\usepackage[]{graphicx}

\renewcommand{\vec}[1]{{\bf #1}}
\newcommand{\scalarproduct}[2]{\left\langle #1 , #2 \right\rangle}

\newcommand{\beq}{\begin{equation}}
\newcommand{\eeq}{\end{equation}}
\newcommand{\barr}{\begin{array}}
\newcommand{\earr}{\end{array}}

\newcommand{\mr}[1]{{\mathrm{#1}}}
\newcommand{\CL}{{C_\mr{L}}}
\newcommand{\CD}{{C_\mr{D}}}

\newcommand{\FL}{{F_\mr{L}}}

\newcommand{\va}{{v_\mr{a}}}
\newcommand{\vasq}{{v^2_{\mr{a}}}}
\newcommand{\vw}{{ v_\mr{w}}} 
\newcommand{\vwcube}{{v^3_{\mr{w}}}}

\newcommand{\dd}{{\mr{d}}}

\begin{document}
\DOIsuffix{theDOIsuffix}
\Volume{VV}
\Issue{I}
\Month{MM}
\Year{YYYY}
\pagespan{1}{}
\Receiveddate{XXXX}
\Reviseddate{XXXX}
\Accepteddate{XXXX}
\Dateposted{XXXX}
\keywords{Quaternions, optimal control, tethered kite, airborne wind energy, kite power, pumping cycle.}
\subjclass[msc2010]{49N90}



\title[Quaternion-based model for airborne wind energy system]{A quaternion-based model for optimal control of the SkySails airborne wind energy system}


\author[M. Erhard]{Michael Erhard\inst{1,2}%
\footnote{Corresponding author,~e-mail:~\textsf{michael.erhard@skysails.de},
     Phone +49\,40\,70299\,203,
     Fax +49\,40\,70299\,222}}
\address[\inst{1}]{SkySails GmbH, Luisenweg 40, 20537 Hamburg, Germany}
\address[\inst{2}]{Systems Control and Optimization Laboratory,
  Department of Microsystems Engineering (IMTEK),
University of Freiburg, Georges-K\"ohler-Allee 102, 79110 Freiburg, Germany}
\address[\inst{3}]{Department ESAT (STADIUS/OPTEC), KU Leuven University, Kasteelpark Arenberg 10, 
   3001 Leuven, Belgium}
\author[G. Horn]{Greg Horn\inst{2,}
}
\author[M. Diehl]{Moritz Diehl\inst{2,3}
}
\begin{abstract}
Airborne wind energy systems are capable of extracting energy from higher wind speeds at higher altitudes. 
The configuration considered in this paper is based on a tethered kite flown in a pumping orbit. 
This pumping cycle generates energy by winching out at high tether
forces and driving a generator while flying figures-of-eight, or lemniscates, as crosswind pattern. 
Then, the tether is reeled in while keeping the kite at a neutral position, thus leaving a net amount of generated energy. 
In order to achieve an economic operation, optimization of pumping cycles is of great interest.

In this paper, first the principles of airborne wind energy will be
briefly revisited.
The first contribution is a singularity-free model for the tethered
kite dynamics in quaternion representation, where the model is derived from first principles.
The second contribution is an optimal control formulation and
numerical results for complete
pumping cycles. Based on the developed model, the setup of the optimal
control problem (OCP) is described in detail along with its numerical
solution based on the direct multiple shooting method in the CasADi
optimization environment.
Optimization results for a pumping cycle consisting of six lemniscates
show that the approach is capable to find an optimal orbit in a few
minutes of computation time. For this optimal orbit, the
power output is increased by a factor of two compared to a
sophisticated initial guess for the considered test scenario.

\end{abstract}
\maketitle                   






\section{Introduction}
The idea of airborne wind energy (AWE) is to generate usable power by
airborne devices. In contrast to towered wind turbines, 
airborne wind energy systems are flying, 
usually connected by a tether to the ground, like
kites or tethered balloons. AWE systems 
exploit the relative velocity between the airmass and the ground and maintain a tension in the tether.
They can be connected to a stationary
ground station, or 
to another moving, but non-flying object, like a land or sea
vehicle. Power is generated in form of a traction force, e.g. to a
moving vehicle, or 
in form of electricity. The three major reasons why people are interested in airborne wind
energy for electricity production are the following:
\begin{itemize}
\item First, like solar, wind power is one of the 
few renewable energy resources
  that is in principle large enough to satisfy all of humanity's
  energy needs.
\item Second, in contrast to ground-based wind turbines, airborne wind energy devices might be able to
  reach higher altitudes, tapping into a
  large and so far unused wind power resource \cite{Archer2009}. The winds in higher
  altitudes are typically stronger and more consistent than those 
close to the ground, both on- and off-shore.
\item Third, and most important, airborne wind energy systems might 
need less material  investment per unit of usable power than most
  other renewable energy sources. This high power-to-mass ratio  promises to make large scale
  deployment of the technology possible at  comparably low costs.
\end{itemize} 

In order to achieve an economic and competitive operation of airborne
wind energy (AWE)
setups, 
optimal control of energy production is the key to success. 
In this paper, a novel formulation of the dynamics of a specific AWE
system is given, and it is shown how it allows one to solve challenging 
optimal control problems for this system.

The paper is organized as follows:  in Section~\ref{crosswind}, 
we introduce the main ideas of
airborne wind energy, following the lines of~\cite{Diehl2013a}, and then
explain in Section~\ref{skysailssystem} the experimental system 
developed by the company SkySails,
the modeling of which is the main subject of this paper. Starting
point is an existing differential equation model of the
system from~\cite{Erhard2012a}, described in
Section~\ref{sec:setup_reference_model}, which suffers from singularities caused by the
choice of coordinate system. The main contribution of this paper is
presented in Section~\ref{sec:quaternion_model}, where a new singularity free model
based on quaternions is developed. It is demonstrated by
numerical simulation that the old and new model 
describe identical system dynamics when away from the singularities, 
while only the new model can be reliably simulated everywhere. It is then
shown in Section~\ref{sec:optimal_control} how the new quaternion-based model can be used to formulate the 
complex periodic optimal control problem (OCP).
In Sect.~\ref{sec:results}, the discrete formulation of the OCP to be solved in  
the dynamic optimization
environment CasADi~\cite{Andersson2012a} is given. Furthermore, resulting complete operational
pumping cycles for the SkySails airborne wind energy system are presented. 
The paper ends with a conclusion in Section~\ref{sec:conclusion}.

\section{Crosswind kite power and pumping cycles}
\label{crosswind}
Every hobby kite pilot or kite surfer knows this observation: As soon as a kite is
flying fast loops in a crosswind direction the tension in the lines increases significantly. The hobby kite pilots have
to compensate the tension strongly with their hands while kite surfers make use of the enormous crosswind
power to achieve high speeds and perform spectacular stunts. The reason for this observation
is that the aerodynamic lift force $\FL$ of an airfoil
increases with the square of the flight velocity, or more exactly,
with the apparent airspeed at the wing, which we denote by $\va$.
More specifically, 
\begin{equation}\label{liftforceformula}
\FL = \frac{1}{2} \rho A \CL \vasq,
\end{equation}
where $\rho$ is the density
of the air, $A$ the airfoil area, and $\CL$ the lift coefficient which
depends on the geometry of the airfoil. 

Thus, if we fly a kite in crosswind direction with a velocity $\va$  that is
ten times faster than the wind speed $\vw$, the tension in the
line will increase by a factor of hundred in comparison to a kite that
is kept at a static position in the sky. The key observation is now
that the high speed of the kite can be maintained by the ambient wind
flow, and that either the high speed itself or the tether tension can
be made useful for harvesting a part of the enormous amount of power that the
moving wing can potentially extract from the wind field. 

 The idea of power generation with tethered airfoils flying fast in 
a crosswind direction was
already in the 1970’s and 1980’s investigated in detail by the
American engineer Miles Loyd \cite{Loyd1980}. 
He was arguably the first to compute the power generation potential of fast
flying tethered wings - a principle that he termed \textit{crosswind kite
  power}. Loyd investigated (and also patented) the following idea: 
an airplane, or kite, is flying on circular trajectories in the sky while being connected to the ground 
with a strong tether. He described two different ways to 
make this highly concentrated form of wind power useful
for human needs, that he termed \textit{lift mode} 
and \textit{drag mode}: while the lift mode uses the tension in the
line to pull a load on the ground, the drag mode uses the high
apparent airspeed to drive a small wind turbine on the wing. 

This paper is concerned with systems operating in lift mode, which has the advantage
over the drag mode that that it does not need high
voltage electrical power transmission via
the tether. The idea is to directly use the strong tether 
tension to unroll the tether from a
drum, and the rotating drum drives an electric generator. As both the
drum and generator can be placed on the ground, we
call this concept also \textit{ground-based generation} or \textit{traction
power generation}. For continuous operation, one has to periodically
retract the tether. One does so by changing the flight pattern to
one that produces much less lifting force. This allows one to reel in the
tether with a much lower energy investment than what was gained in
the power production phase. The power production phase is also called
\textit{reel-out} phase, and the retraction phase \textit{reel-in} 
or \textit{return} 
phase. Loyd
coined the term \textit{lift
  mode} because one uses mainly the lifting force of the wing. But due
to the periodic reel-in and reel-out motion of the tether, this way of 
ground-based power generation is often also called \textit{pumping
mode}; sometimes even the term \textit{Yo-Yo mode} was used to describe it.

\begin{figure}
\begin{minipage}[b]{50mm}
\centering
\includegraphics[height=7cm]{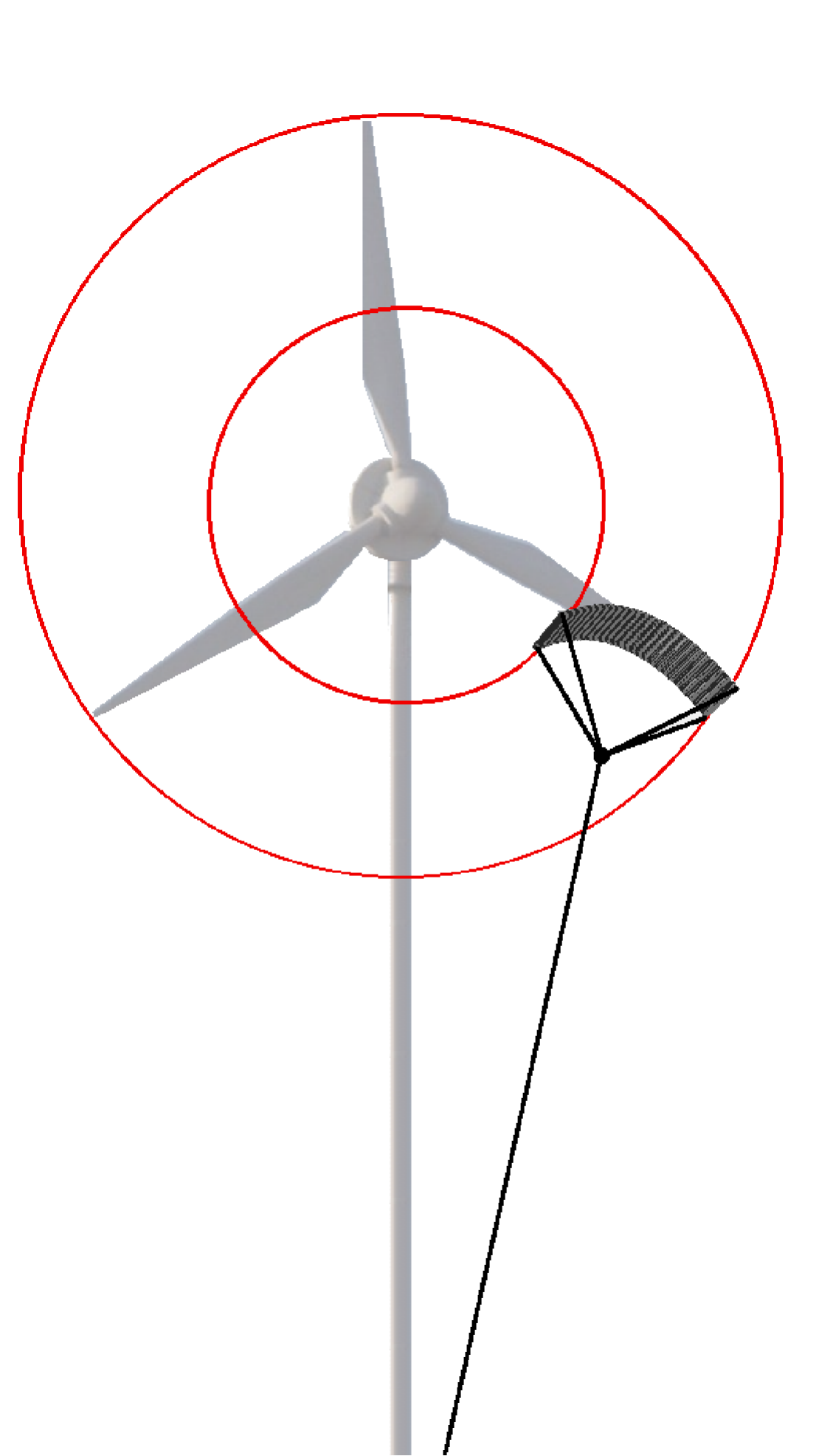}
\caption{The basic idea of the power generating device is to only build the wing
tips of a gigantic windmill, in form of tethered airfoils, i.e.
kites. Fast crosswind motion of the kites creates an enormous
traction force
[graphics courtesy Boris Houska].}
\label{windmil3}
\end{minipage}
\hfil
\begin{minipage}[b]{90mm}
\includegraphics[width=\linewidth]{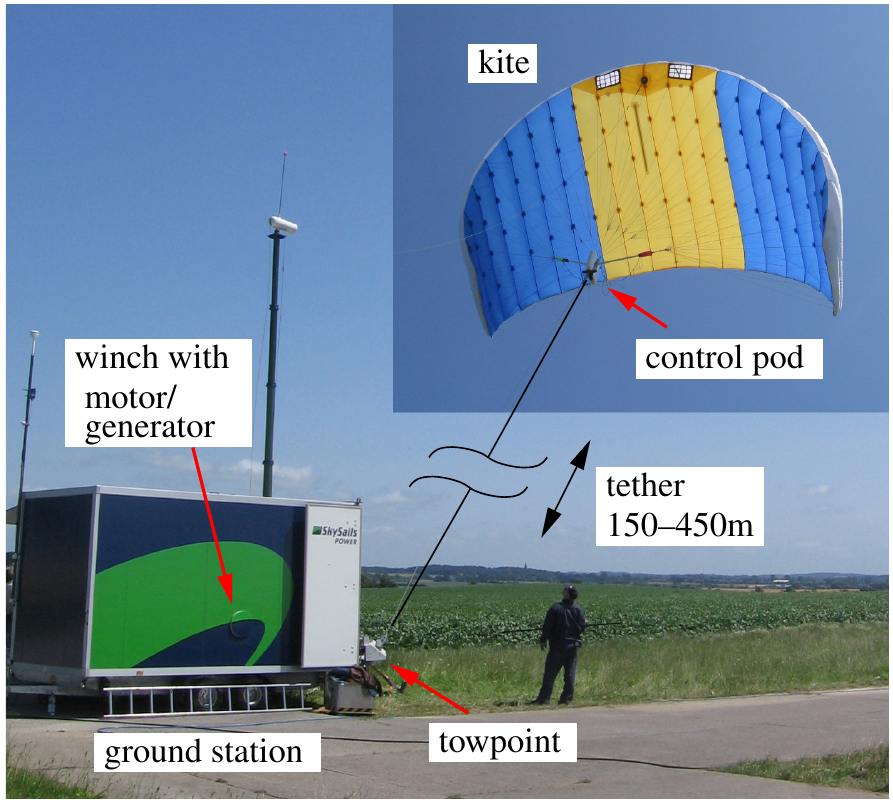}%
\caption{SkySails functional prototype setup for kites of sizes ranging from 20 to 40\,m$^2$ (30\,m$^2$ shown here). A tether in the range 150--400\,m transfers the forces from the flying system to the ground based main winch, which is attached to a motor/generator. A specific feature of this setup is the single tether and the airborne control pod located under the kite. The steering actuator in the control pod pulls certain lines in order to steer the kite.}
\label{fig:sks_setup}
\end{minipage}
\end{figure}
It is interesting to compare crosswind kite power systems with a
conventional wind turbine, as done in
Fig.~\ref{windmil3}, which shows a conventional
wind turbine superposed with an airborne wind energy
system.
Seen from this perspective, the idea of AWE is to only
build the fastest moving part of a large wind turbine, the tips of
the rotor blades, and to replace them by automatically controlled fast flying
tethered kites. The motivation for this is the fact that the outer
30\% of the blades of a conventional wind turbine provide more than
half of the
total power, while they are much thinner and lighter than the inner parts of
the blades. Roughly speaking, the idea of airborne wind energy systems
is to replace the inner parts of the blades, as well as the
tower, by a tether. 

The power $P$ that can be generated with a tethered airfoil operated
either in drag or in lift mode had under idealized assumptions been
estimated by Loyd \cite{Loyd1980} to
be approximately given by 
\begin{equation}\label{loydsformula}
P = \frac{2}{27} \rho A \vwcube \CL \left(\frac{\CL}{\CD}\right)^2,
\end{equation}
where $A$ is the area of the wing, $\CL$ the lift and $\CD$
the drag coefficients, and $\vw$ the wind speed. Note that the
lift-to-drag ratio $E \doteq \frac{\CL}{\CD}$ enters the formula quadratically and is
thus an important wing property for crosswind AWE systems.  For
airplanes, this ratio is also referred to as the \textit{gliding
  number}; it describes how much faster a glider without propulsion can move horizontally
compared to its vertical sink rate.

Since the first proposal of using tethered wings for energy harvesting
by Loyd,
academic research and industrial development activities in the field
of airborne wind energy (AWE) have been 
significantly increasing, especially during the last 10 years.
The main advantage of the AWE technology is that
airborne systems are 
capable of harvesting energy from higher wind speeds at higher altitudes. 
In addition, as most setups require less installation efforts compared
to conventional wind turbines, e.g.~foundations and towers, 
this technology is a promising candidate to become an additional
source of 
renewable energy in the near future.
For an overview on the different geometries which have been developed
so far, the reader is referred to \cite{Fagiano2012a} and 
the monograph on AWE \cite{Ahrens2013}.

\section{The SkySails Power prototype}
\label{skysailssystem}
The geometry under consideration in this paper is based on the SkySails Power prototype, shown in Fig.~\ref{fig:sks_setup}.
This setup features an airborne control pod located under the soft kite connected to a ground based winch with motor/generator by a single tether.
The functional prototype successfully demonstrated fully automated energy generation by autonomously controlled pumping cycles \cite{Erhard2015a}.
\begin{vchfigure}
\includegraphics[width=7cm]{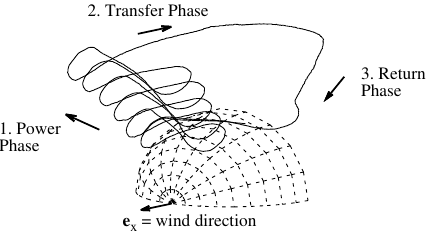}
\vchcaption{Trajectory of an experimentally flown pumping cycle \cite{Erhard2015a}.}
\label{fig:pumping_cycle}
\end{vchfigure}
The principle of energy generation is depicted in Fig.~\ref{fig:pumping_cycle} and shall be briefly summarized:
during the power phase, the kite is flown dynamically in lemniscates
(figures of eight). This dynamical crosswind flight leads to high tether forces. By reeling out the tether, electrical energy is produced.
At a certain line length, the kite is transferred to the neutral zenith position. 
At this low-force position, the tether is reeled in (return phase) consuming a certain portion of the previously generated energy. 
Finally, a reasonable amount of generated net power remains. Hence,
the average power output over complete cycles will be the target of the optimization implementation presented in this paper.

The optimal control is based on a dynamical model of the system, and
the chosen numerical techniques are similar to the ones used in~\cite{Houska2007}. In
order to allow for optimization of complete pumping cycles, a simple
model should be chosen in order to keep the number of optimization variables low.
The model developed in this paper is based on a very simple model set
up for the 
SkySails system \cite{Erhard2012a}, which describes the basic dynamic properties of the system and has been experimentally validated \cite{Erhard2013b}.
The quaternion representation allows for singularity-free equations of motion, and the replacement of trigonometric functions seems to be advantageous for the optimization process.

A further advantage of the considered system is that an experimentally flown trajectory (see e.g.~Fig.~\ref{fig:pumping_cycle}) can be used as initial guess and the optimization result can be directly compared to it in return.
It should be noted that the optimization has to be 
subject to certain geometrical constraints, which will be discussed in this paper.
Particularly, topological constraints have to be imposed to preserve the 
topology of the pumping cycle during the optimization process (the 6 lemniscates given in the initial guess, compare Fig.~\ref{fig:pumping_cycle}).
\section{Setup and reference kite model}
\label{sec:setup_reference_model}
\begin{figure}
  \centering
  \sidecaption
  \includegraphics[width=5cm]{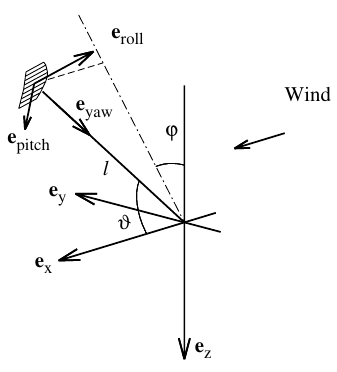}
  \caption{Coordinate system definitions of the fixed reference vectors $\vec{e}_x, \vec{e}_y, \vec{e}_z$ and the body frame vectors 
  $\vec{e}_{\rm roll}$, $\vec{e}_{\rm pitch}$ and   $\vec{e}_{\rm yaw}$.
  The position is determined by the two angles $\varphi$, $\vartheta$ and the
  tether length $l$.
  The wind direction is defined in $\vec{e}_x$-direction.}
  \label{fig:coords_definition}
\end{figure}

\begin{figure}
\begin{minipage}{75mm}
\includegraphics[width=\linewidth]{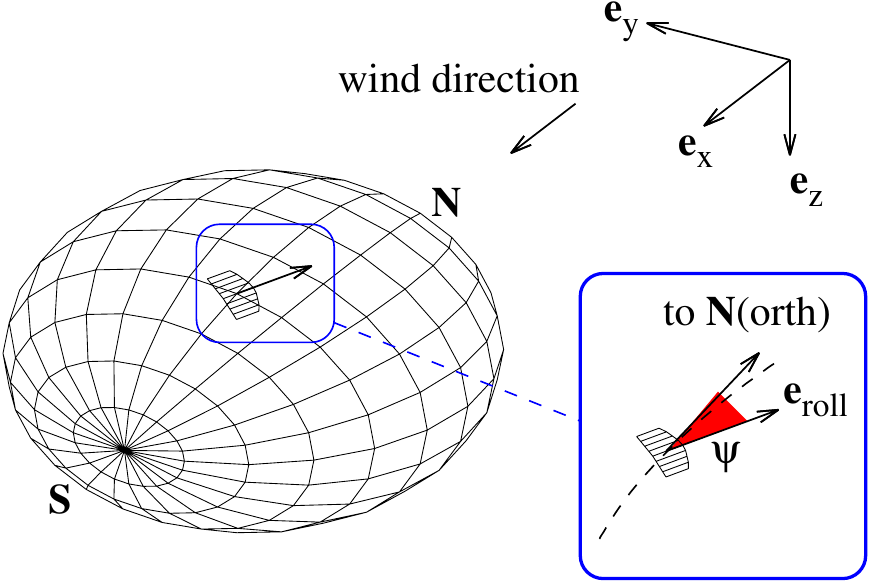}%
\caption{Demonstrative explanation of the used Euler angles by
  comparing them to angles used for navigation on earth. For that
  purpose, the earth's rotation axis has to be tilted by 90 degrees
  with the north pole towards the wind direction. The angles
  $\varphi$, $\vartheta$ 
correspond to the position on earth and $\psi$ to the bearing.
The exact mapping is given in Table \ref{tab:coords}.}
\label{fig:coords}
\end{minipage}
\hfil
\begin{minipage}{80mm}
\setfloattype{table}
\caption{Mapping of the Euler angles of the tethered system to earth navigation quantities (compare Fig.~\ref{fig:coords}).}
\label{tab:coords}
\begin{tabular}{cl}
  \hline\hline
 angle in model &  'earth' navigation quantity \\
  \hline
$\varphi$ &   earth longitude \\
$\vartheta$ &    earth latitude \\
   $\psi$ & bearing (direction w.r.t.~north pole)\\
  \hline\hline
\end{tabular}
\end{minipage}
\end{figure}

In this section, the model from \cite{Erhard2012a}, \cite{Erhard2013b} is briefly summarized.
The coordinate system is defined in Fig.~\ref{fig:coords_definition}.
The state vector of the dynamical system consists of the 3-dimensional orientation of the
flying system, given by $\left\{\psi, \varphi, \vartheta \right\}$ and the
tether length $l$.
\begin{equation}
  {\bf x} = \left[ \psi, \varphi, \vartheta, l \right]^\top 
\end{equation}
Note that the tether of the system leads to a relation between orientation and
position, which is given in the basis $\{\vec{e}_x,\vec{e}_y,\vec{e}_z\}$ by:
\begin{equation}
  {\bf r} = l\left[
    \begin{array}{c}
      \cos\vartheta \\ \sin\varphi \sin\vartheta \\ -\cos\varphi \sin\vartheta
    \end{array} 
  \right]\label{eq:r_from_angles}
\end{equation}
The wind direction is defined in $\vec{e}_x$-direction.
In order to provide a clear picture of the coordinate system definition, a comparison with the navigation task on earth shall be discussed as illustrated in Fig.~\ref{fig:coords}.
The coordinate system corresponds to an earth coordinate system with
symmetry (rotational) axis in wind direction, which could be achieved
by rotating the north pole by 90 degrees towards the wind direction. 
According to the angle mapping given in Table \ref{tab:coords}, the two angles 
$\varphi, \vartheta$ correspond to longitude and latitude, hence determine the position on earth. 
The orientation angle $\psi$ corresponds to bearing, i.e.~direction w.r.t.~north in our earth system. This angle determines the direction of motion (course) and
the angle of $\psi=0$ describes heading directly against the wind.
It should be finally noted that due to the ambient wind, the angle $\psi$ does not exactly coincide the kinematic flight direction (course). This effect is fully covered by the developed model and further explained in \cite{Erhard2015a}.

The input vector of the model comprises the inputs of the two control actuators: the
steering deflection of the airborne control pod, $\delta$ and the winch speed of
the tether $v^{\rm (winch)}$, respectively.
\begin{equation}
  {\bf u} = \left[ \delta, v^{\rm (winch)} \right]^\top
\end{equation}

In the following, the equations of motion shall be summarized for
reference purpose only. This model was derived based on the same
assumptions as those that will be used in the subsequent section for the derivation of the quaternion representation.
As conducted in detail in \cite{Erhard2012a}, \cite{Erhard2013b}, the equations of motion read:
\begin{eqnarray}
  \dot{\psi} &=& g_{\rm k}\, v_{\rm a}\, \delta + \dot{\varphi}
  \cos\vartheta \label{eq:eqm_psi}\\
  \dot{\varphi} &=& -\frac{v_{\rm a}}
  {l \sin\vartheta} \sin\psi\label{eq:eqm_phi}\\
  \dot{\vartheta} &=& 
  - \frac{v_{\rm w}}{l}
  \sin\vartheta
  +\frac{v_{\rm a}}{l}\cos\psi \label{eq:eqm_theta}\\
  \dot{l} &=& v^{\rm (winch)}\label{eq:eqm_l} 
\end{eqnarray}
Here, the airpath speed of the kite, $v_{\rm a}$, is given by
\begin{equation}
  v_{\rm a} = v_{\rm w}E \cos\vartheta - \dot{l}E . \label{eq:eqm_v_a}
\end{equation}
The set of equations involves three parameters: lift-to-drag (glide) ratio $E$, steering response proportionality constant $g_{\rm k}$ and ambient wind speed $v_{\rm w}$.
Further explanations of parameters including values used later in the paper are summarized in Table~\ref{tab:parameters}.
\begin{vchtable}
  \vchcaption{Description of system parameters and values used for optimization}
  \label{tab:parameters}
  \begin{tabular}{@{}lllp{8cm}@{}}
    \hline
    symbol & \multicolumn{2}{l}{value} & description \\
    \hline
    $A$ & 21.0 & m$^2$ & projected area of the kite\\
    $C_{\rm R}$ & 1.0 && aerodynamic force coefficient, given by
    $C_{\rm D}=\sqrt{C_{\rm L}^2+C_{\rm D}^2}$ for aerodynamic lift (drag) coefficient
    $C_{\rm L}$ ($C_{\rm D}$)\\
    $\dot{\delta}_{\rm max}$ & 0.6 & 1/s & maximal control actuator steering speed \\
    $\delta_{\rm max}$ & 0.7 && maximal control actuator steering deflection\\
    $E$ & 5.0 && lift-to-drag (glide) ratio. Note, that in contrast to most other systems, the SkySails system is operated at constant glide ratio\\
    $g_{\rm k}$& 0.1 & rad/m & proportionality constant relating turn rate of
    the kite to steering actuator deflection\\
    $l_{\rm max}$ & 300 & m & maximal available tether length\\
    $\theta_{\rm min}$ & 0.35 & rad & minimal elevation angle w.r.t.~the horizontal plane\\
    $\varrho$ & 1.2 & kg/m$^3$& air density\\
    $v_{\rm a,min}$ & 5.0 & m/s & minimal air path speed in order to guarantee
    flight stability and keep system tethered (this corresponds to free flight velocity)\\
    $v^{\rm (winch)}_{\rm min}$ & $-5.0$ & m/s & typically chosen as $-v_{\rm w}/2$, leads to a windward limit of $\approx 110^\circ$\\
    $v_{\rm w}$ & 10 & m/s & ambient wind speed\\
    \hline
  \end{tabular}
 \end{vchtable}
Finally, the tether force can be computed based on the air path speed $v_{\rm a}$ by 
\begin{equation}
  F_{\rm tether} = \frac{\varrho A C_{\rm R}}{2}\frac{E}{\sqrt{1+E^2}} v_{\rm a}^2
\end{equation}
with density of air $\varrho$, projected kite area $A$ and aerodynamic force coefficient $C_{\rm R}$.
\section{Quaternion based kite model}
\label{sec:quaternion_model}
In the following, the equations of motion shall be derived
in a singularity-free formulation based on quaternions.
Quaternions are a well-known tool in aerospace modeling and
estimation to avoid the singularities of Euler
angles~\cite{Kuipers1999}.
An alternative singularity-free representation could be based on
rotation matrices 
as natural coordinates directly, as is shown for modeling of AWE systems in~\cite{Gros2014a}. Note that the rotation matrix in the following will be used only as intermediate step in the derivation towards quaternions.
The base coordinate system is depicted in Fig.~\ref{fig:coords}a.
The time dependence w.r.t.~to this reference frame is expressed by rotation matrices $R(t) \in \mathrm{SO}(3)$, i.e.
\begin{eqnarray}
  \vec{e}_{\rm roll}(t) &=& -R(t) \vec{e}_{z} \label{eq:def_e_roll}\\
  \vec{e}_{\rm pitch}(t) &=& -R(t) \vec{e}_{y} \\
  \vec{e}_{\rm yaw}(t) &=& -R(t) \vec{e}_{x} \label{eq:def_e_yaw}
\end{eqnarray}

\begin{figure}
  \sidecaption
  \includegraphics[width=8cm]{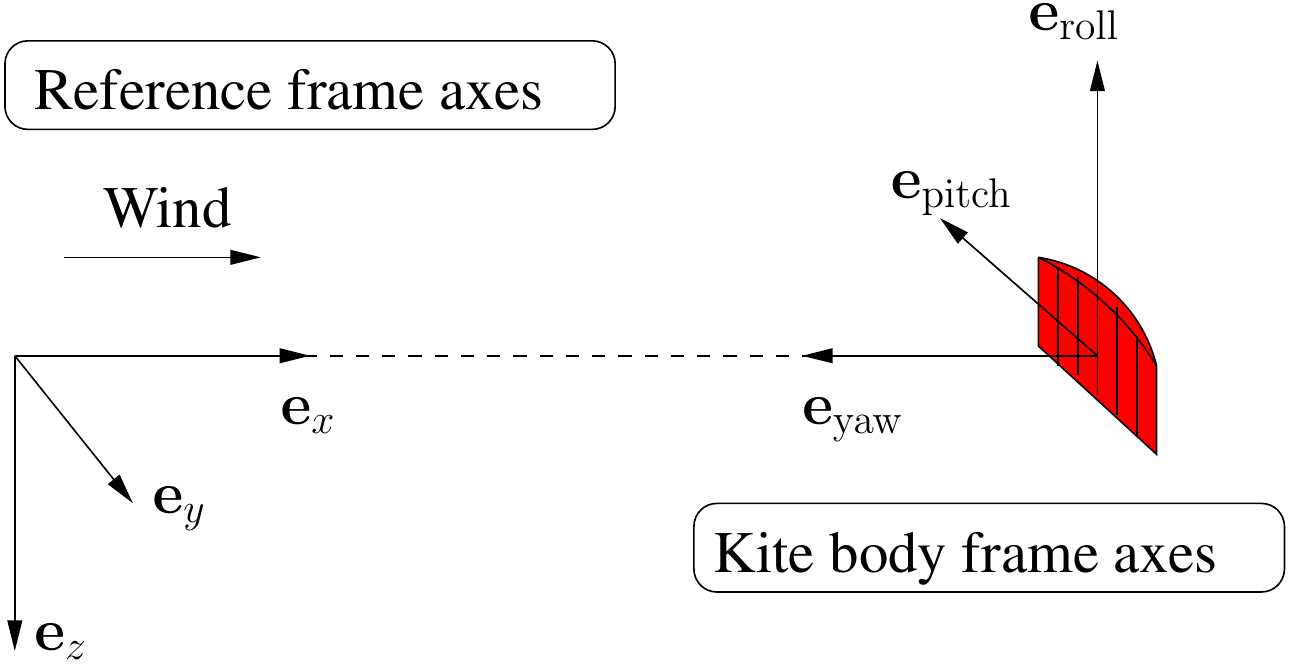}
  \caption{Reference frame axes for the derivation of the quaternion equations of motion. The depicted orientation corresponds to the zero-rotation, 
  i.e.~$R=I$ or $\vec{q}=[1,0,0,0]$.}
  \label{fig:zero_rotation}
\end{figure}

The zero-rotation state is defined by 
$R \doteq I$, corresponding to 
$\vec{e}_{\rm yaw}=-\vec{e}_x$,
$\vec{e}_{\rm pitch}=-\vec{e}_y$ and 
$\vec{e}_{\rm roll}=-\vec{e}_z$ as illustrated in 
Fig.~\ref{fig:zero_rotation}. 
Note that this seems to be an arbitrary choice 
a priori. However defining the ambient wind direction, given by $\vec{e}_x$, 
in the kite's yaw-axis leads to a beneficial simplicity of the equations
of motion later.

It should be emphasized that due to the nature of the tethering, position
and orientation are related. Hence, the position of the kite is given by
\begin{equation}
  \vec{r}(t) = -l(t) \vec{e}_{\rm yaw}(t) = l(t) R(t) \vec{e}_{x}
  \label{eq:def_position}
\end{equation}
The air path speed vector can be calculated as sum of the ambient wind
vector $v_{\rm w}\vec{e}_x$ and reversed kinematic velocity
vector $-\dot{\vec{r}}$ (apparent wind opposed to motion). Using the time
derivative of (\ref{eq:def_position}) and introducing the speeds $v_{\rm roll}$ and $v_{\rm pitch}$ in the tangent plane, compare Fig.~\ref{fig:tangent_plane},  yields 
\begin{equation}
  \vec{v}_{\rm a}(t) = 
    v_{\rm w} \vec{e}_x - 
    v_{\rm pitch}(t) \vec{e}_{\rm pitch}(t) -
    v_{\rm roll}(t) \vec{e}_{\rm roll}(t) +
    \dot{l}(t) \vec{e}_{\rm yaw}(t) \label{eq:def_vector_v_a}
\end{equation}

\subsection{Model assumptions}
For AWE systems, the aerodynamic forces are typically large compared to masses, which in fact is a prerequisite for tethered flight. This is particularly true for the SkySails Power system.
Therefore, acceleration effects play a minor role and can be neglected. In addition, the system is assumed to fly always in its aerodynamic equilibrium state. These two simplifications lead to the simple structure of the equations of motion.
\subsubsection{Aerodynamics}
\begin{figure}
  \sidecaption
  \includegraphics[width=7cm]{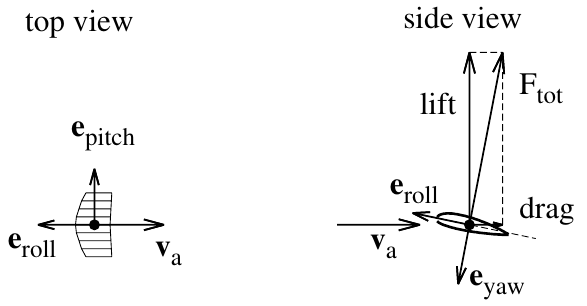}
  \caption{The kite is assumed to always fly in an aerodynamic equilibrium
  state. As aerodynamic model, absence of transversal air flow (top view) and     constant lift-to-drag-ratio (side view) is  demanded. Figure
  taken from \cite{Erhard2013b}}
  \label{fig:aerodynamics}
\end{figure}
The aerodynamics of the system is reduced to two conditions, depicted in
Fig.~\ref{fig:aerodynamics} and given in the following.
\begin{enumerate}
  \item The kite experiences no transversal wind flow, i.e.~no 'side-slip' angle occurs, compare Fig.~\ref{fig:aerodynamics}{\bf b}.
  \begin{equation}
    \scalarproduct{\vec{e}_{\rm pitch}(t)}{\vec{v}_{\rm a}(t)} = 0
  \end{equation}
  Inserting (\ref{eq:def_vector_v_a}) yields
  \begin{equation}
    \scalarproduct{\vec{e}_{\rm pitch}(t)}{v_{\rm w}\vec{e}_x} - 
    v_{\rm pitch}(t) +
    \dot{l} \underbrace{
    \scalarproduct{\vec{e}_{\rm pitch}(t)}{\vec{e}_{\rm yaw}(t)}}_{=0} = 0
  \end{equation}
  Thus, the condition results in
  \begin{equation}
    v_{\rm pitch}(t) = 
    -v_{\rm w} \scalarproduct{R(t)\vec{e}_{y}}{\vec{e}_x}
  \end{equation}
  \item The lift-to-drag ratio is constant and given by the parameter $E$. Hence, the kite is assumed to fly always at the same angle of attack. Considering the geometry of Fig.~\ref{fig:aerodynamics}{\bf c}) yields
  \begin{equation}
    \scalarproduct{\vec{e}_{\rm roll}(t)}{\vec{v}_{\rm a}(t)} = E
    \scalarproduct{\vec{e}_{\rm yaw}(t)}{\vec{v}_{\rm a}(t)}
  \end{equation}
  Inserting (\ref{eq:def_vector_v_a}) produces
  \begin{equation}
    v_{\rm roll}(t) =
    v_{\rm w}\scalarproduct{R(t)(E \vec{e}_{x} - \vec{e}_{z})}{\vec{e}_x}
    - E \dot{l}(t)
  \end{equation}
  \end{enumerate}
Finally, by considering the motion in the tangent plane, compare Fig.~\ref{fig:tangent_plane}, the following relations for the
turn rates in kite-body-frame can be derived.
\begin{vchfigure}
  \sidecaption
  \includegraphics[width=3.5cm]{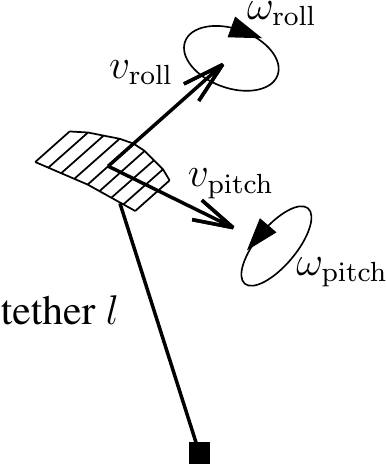}
  \caption{Motion in the tangent plane. Due to the tether, the velocities in the tangent plane $v_{\rm roll}$ and $v_{\rm pitch}$ determine $\omega_{\rm roll}$ and $\omega_{\rm pitch}$ as given in (\ref{eq:omega_roll}) and (\ref{eq:omega_pitch}).}
  \label{fig:tangent_plane}
\end{vchfigure}
\begin{eqnarray}
  \omega_{\rm roll} &=& \frac{v_{\rm pitch}(t)}{l(t)} 
  = -\frac{v_{\rm w}}{l(t)} 
  \scalarproduct{R(t)\vec{e}_{y}}{\vec{e}_x}
  \label{eq:omega_roll}\\
  \omega_{\rm pitch} &=& -\frac{v_{\rm roll}(t)}{l(t)} 
  = - \frac{v_{\rm w}}{l(t)} 
  \scalarproduct{R(t)(E \vec{e}_{x} - \vec{e}_{z})}{\vec{e}_x}
    + E\frac{\dot{l}(t)}{l(t)}
  \label{eq:omega_pitch}
\end{eqnarray}
\subsubsection{Steering}
Steering of the kite is accomplished by the control actuator in the control pod.
Applying an actuator deflection $\delta$ to the steering lines tilts the kite canopy resulting in a curve flight or turn rate around the $\vec{e}_{\rm yaw}$-axis.
The steering behavior can be described by the following equation
\begin{equation}
  \omega_{\rm yaw} = g_{\rm k} v_{\rm a}\delta \label{eq:omega_yaw}
\end{equation} 
This \textit{turn-rate law} states that the rotation rate is proportional to a constant system parameter $g_{\rm k}$, the deflection $\delta$ and the scalar air path speed $v_{\rm a}$, which is defined as component of $\vec{v}_{\rm a}$ in $-\vec{e}_{\rm roll}$-direction
\begin{equation}
  v_{\rm a} \doteq 
  -\scalarproduct{\vec{e}_{\rm roll}(t)}{\vec{v}_{\rm a}(t)} =
  \scalarproduct{R(t) \vec{e}_{z}}{\vec{v}_{\rm a}(t)}
\end{equation}
The validity of the turn-rate law has been experimentally shown for different kites \cite{Erhard2012a}, \cite{Fagiano2013a}, \cite{Jehle2012}. In addition, a nice derivation based on geometric considerations can be found in \cite{Fagiano2013a}.
It should be mentioned that correction terms for mass effects \cite{Erhard2013a} and aerodynamic effects \cite{Costello2013a} might be added, but are neglected in this paper.
\subsection{Quaternions}
The final formulation of the equations of motion is based on quaternions~\cite{Kuipers1999}, introduced as
\begin{equation}
  \vec{q} = \left[\begin{array}{c}
    q_0 \\ q_1 \\ q_2 \\ q_3 
  \end{array}\right] 
\end{equation} 
The time evolution as function of turn rates, which are given in the body frame, reads 
\begin{equation}
  \dot{\vec{q}} = \frac{1}{2} \left[
  \begin{array}{cccc}
    0 & -\omega_{x} & -\omega_{y} & -\omega_{z} \\
    \omega_{x} & 0 &   \omega_{z} & -\omega_{y} \\
    \omega_{y} & -\omega_{z} & 0 & \omega_{x} \\
    \omega_{z} & \omega_{y} & -\omega_{x} & 0
  \end{array}
  \right] \vec{q}
\end{equation}
Mapping of the kite axes to $\omega_x,\omega_y,\omega_z$ has to be considered for $R\!=\!I$. 
Regarding (\ref{eq:def_e_roll}--\ref{eq:def_e_yaw}) or
Fig.~\ref{fig:zero_rotation} yields  
$\omega_x\!=\!-\omega_{\rm yaw}$,
$\omega_y\!=\!-\omega_{\rm pitch}$ and
$\omega_z\!=\!-\omega_{\rm roll}$.
The quaternion time evolution then results in
\begin{equation}
  \dot{\vec{q}} = \frac{1}{2} \left[
  \begin{array}{cccc}
    0 & \omega_{\rm yaw} & \omega_{\rm pitch} & \omega_{\rm roll} \\
    -\omega_{\rm yaw} & 0 &   -\omega_{\rm roll} & \omega_{\rm pitch} \\
    -\omega_{\rm pitch} & \omega_{\rm roll} & 0 & -\omega_{\rm yaw} \\
    -\omega_{\rm roll} & -\omega_{\rm pitch} & \omega_{\rm yaw} & 0
  \end{array}
  \right] \vec{q}\label{eq:eqm_quaternions}
\end{equation}
The relation to the rotation can be given by
\begin{equation}
  R(t) = \left[
  \begin{array}{ccc}
    q_0^2+q_1^2-q_2^2-q_3^2 &
    2 (q_1 q_2-q_0 q_3) &
    2 (q_1 q_3 + q_0 q_2) \\[8pt]
    2 (q_1 q_2 + q_0 q_3) &
	q_0^2 - q_1^2 + q_2^2 - q_3^2 &
    2 (q_2 q_3 - q_0 q_1) \\[8pt]
    2 (q_1 q_3 - q_0 q_2) &
    2 (q_2 q_3 + q_0 q_1) &
    q_0^2 - q_1^2 - q_2^2 + q_3^2
  \end{array}
  \right]\label{eq:R_quaternion}
\end{equation}
The derivation of the system's equations of motion is done as follows:
The turn rates in kite-body-frame,
$\omega_{\rm roll}$, $\omega_{\rm pitch}$ and $\omega_{\rm yaw}$ are 
expressed as functions of $\vec{q}$ and subsequently inserted in (\ref{eq:eqm_quaternions}).
In order to achieve this, (\ref{eq:R_quaternion}) is used for
$R(t)$ and inserted into (\ref{eq:omega_roll}), (\ref{eq:omega_pitch}) and
(\ref{eq:omega_yaw}).
\subsection{Equations of motion}
The state vector for the quaternion formulation is given by 
\begin{equation}
\left[q_0,q_1,q_2,q_3,l\right]^\top
\end{equation}
and the input vector by
\begin{equation}
  \left[\delta, v^{\rm (winch)}\right]^\top
\end{equation}
Performing the derivation steps of the previous section yields the following set of equations of motion:
\begin{eqnarray}
  \dot{\vec{q}} &=& \frac{v_{\rm a}}{2 l} \left[\begin{array}{c}
  -q_2 \\ -q_3 \\ q_0 \\ q_1
  \end{array}\right]
  + \frac{v_{\rm w}}{l} \left[ \begin{array}{c}
    q_0 (q_2^2 + q_3^2) \\
    q_1 (q_2^2 + q_3^2) \\
  - q_2 (q_0^2 + q_1^2) \\
  - q_3 (q_0^2 + q_1^2)
  \end{array}\right] 
  + \frac{g_{\rm k} v_{\rm a} \delta}{2} \left[
  \begin{array}{c}
    q_1 \\ -q_0 \\ -q_3 \\ q_2
  \end{array}
  \right] \label{eq:q_eqm}
  \\
  \dot{l} &=& v^{\rm (winch)}
\end{eqnarray}
with air path speed
\begin{equation}
  v_{\rm a} = E v_{\rm w} (q_0^2+q_1^2-q_2^2-q_3^2) - E \dot{l} \label{eq:q_va}
\end{equation}

It should be remarked that if the winch speed is defined relative to the wind speed by substituting $\dot{l} = v_{\rm w} v^{\rm winch}_{\rm rel}$, the product of 
wind speed $v_{\rm w}$ and time $t$ (or cycle time $T$) becomes an invariant.
As a consequence, solutions for different wind speeds are simply obtained by scaling the dynamic's time axis proportional to $1/v_{\rm w}$.

Although these equations of motion (EQM) preserve the 
quaternion norm $(\mathrm{d}/\mathrm{d}t)\|\vec{q}(t)\|=0$, a stabilization term
in order to compensate for numerical inaccuracies should be added.
The EQM for the numerical optimization read then
\begin{equation}
   \dot{\vec{q}} = \Phi(\cdot) - \gamma_q (\|\vec{q}\|^2-1) \vec{q}
\end{equation}
where $\Phi(\cdot)$ denotes the equations of motion (\ref{eq:q_eqm}--\ref{eq:q_va}) and the damping constant $\gamma_q$ is chosen to lead to a slow decay ($\gamma_q=0.01$\,1/s for the parameters given in Table \ref{tab:parameters}).
\subsection{Transformation to Euler angles for comparison}
For sake of completeness and in order to allow for a simple geometric interpretation of the results in the Cartesian coordinate system (see Fig.~\ref{fig:coords}),
transformations from the set of Euler angles $\{\varphi,\vartheta,\psi\}$ to quaternions and back shall be summarized.
The sequence of rotations
\begin{equation}
  R = R_x(\varphi) R_y(\vartheta) R_x(-\psi)
\end{equation}
is considered in quaternion representation 
\begin{eqnarray}
  \vec{q} &=& \left[
  \begin{array}{c}
    \cos\left(\frac{\varphi}{2}\right)\\
    \sin\left(\frac{\varphi}{2}\right)\\
    0\\
    0    
  \end{array}\right] \otimes
  \left[\begin{array}{c}
    \cos\left(\frac{\vartheta}{2}\right)\\
    0\\
    \sin\left(\frac{\vartheta}{2}\right)\\
    0
  \end{array}\right] \otimes
  \left[\begin{array}{c}
    \cos\left(\frac{\psi}{2}\right)\\
    -\sin\left(\frac{\psi}{2}\right)\\
    0\\
    0
  \end{array}\right]
  \end{eqnarray}
Calculating the quaternion product ($\otimes$) \cite{Kuipers1999}
of the single rotations yields 
\begin{eqnarray}
  \vec{q} &=& \left[
  \begin{array}{c}
  \cos\left(\frac{\varphi}{2}\right)
  \cos\left(\frac{\vartheta}{2}\right)
  \cos\left(\frac{\psi}{2}\right) 
  +\sin\left(\frac{\varphi}{2}\right)
  \cos\left(\frac{\vartheta}{2}\right)
  \sin\left(\frac{\psi}{2}\right) \\
  \sin\left(\frac{\varphi}{2}\right)
  \cos\left(\frac{\vartheta}{2}\right)
  \cos\left(\frac{\psi}{2}\right) 
  -\cos\left(\frac{\varphi}{2}\right)
  \cos\left(\frac{\vartheta}{2}\right)
  \sin\left(\frac{\psi}{2}\right) \\
  -\sin\left(\frac{\varphi}{2}\right)
  \sin\left(\frac{\vartheta}{2}\right)
  \sin\left(\frac{\psi}{2}\right) 
  +\cos\left(\frac{\varphi}{2}\right)
  \sin\left(\frac{\vartheta}{2}\right)
  \cos\left(\frac{\psi}{2}\right) \\
  \sin\left(\frac{\varphi}{2}\right)
  \sin\left(\frac{\vartheta}{2}\right)
  \cos\left(\frac{\psi}{2}\right) 
  +\cos\left(\frac{\varphi}{2}\right)
  \sin\left(\frac{\vartheta}{2}\right)
  \sin\left(\frac{\psi}{2}\right) 
  \end{array}\right]
\end{eqnarray}
The kite position  w.r.t.~$\{\vec{e}_x, \vec{e}_y, \vec{e}_z\}$ can be found by inserting (\ref{eq:R_quaternion}) into (\ref{eq:def_position}) as
\begin{equation}
  \vec{r} = l \left[ \begin{array}{c}
  (q_0^2+q_1^2-q_2^2-q_3^2) \\
  2 (q_0 q_3 + q_1 q_2) \\ 
  2 (q_1 q_3 - q_0 q_2)
  \end{array}\right]\label{eq:r_from_q}
\end{equation}
The Euler angles can be computed from the quaternions by
\begin{eqnarray}
  \varphi &=& \mathrm{atan2}\left((q_0 q_3 + q_1 q_2) , (q_0 q_2 - q_1 q_3)\right) \label{eq:q2varphi} \\
  \vartheta &=& \arccos(q_0^2+q_1^2-q_2^2-q_3^2)\\
  \psi &=& \mathrm{atan2}\left((q_0 q_3 - q_1 q_2) , (q_0 q_2 + q_1 q_3)
  \right)\label{eq:q2position}
\end{eqnarray}
where $\mathrm{atan2}$ is the four-quadrant arctangent function. 
\subsection{Simulation example}
In order to complete the presentation of the quaternion-based model, a simulation example will be discussed. The simulations have been carried out using the quaternion-based model on the one and the reference model from Sect.~\ref{sec:setup_reference_model} on the other hand.
Both computations have employed a Runge-Kutta-4 ODE propagation using timesteps of $\tau\!=\!0.1$\,s, a tether length of $l\!=\!100$\,m and the initial conditions $\varphi(0)\!=\!0, \vartheta(0)\!=\!\arctan(E)$ and $\psi(0)=0$. A constant steering command of $\delta\!=\!0.186$ has been used to fly next to the singularity of the Euler angle representation.
\begin{vchfigure}
\includegraphics[width=12cm]{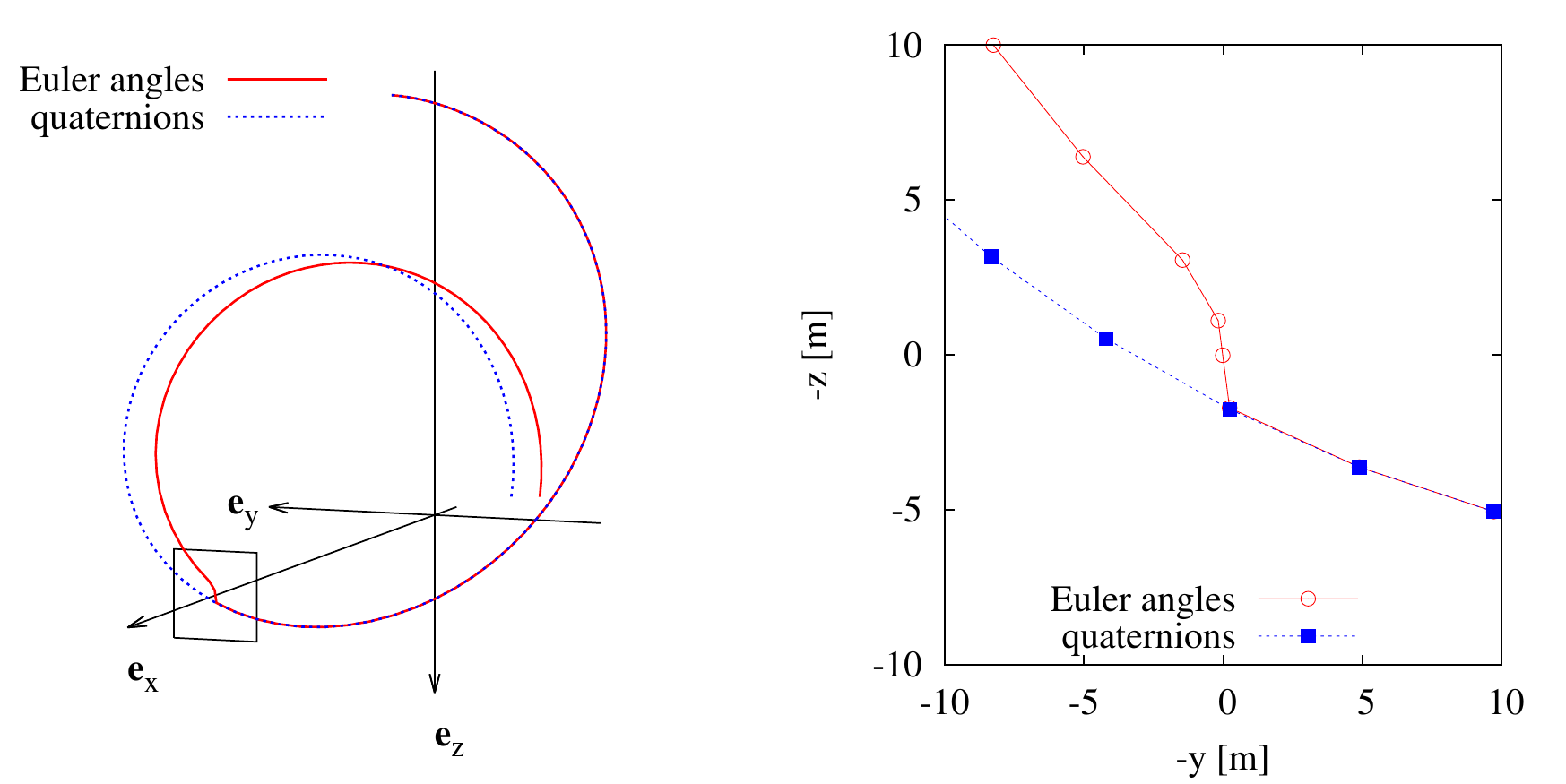}%
\vchcaption{Simulation of the simple model ODE, formulated in Euler angles (solid red) and quaternions (dashed blue). The 3d trajectory is shown on the left, the right plot zooms into the vicinity of the singularity, where the Euler angle formulation fails. The points in the right subfigure represent the discrete simulation timesteps of constant duration.}
\label{fig:sim_singularity}
\end{vchfigure}
The simulation results are depicted in Fig.~\ref{fig:sim_singularity}. 
It can be recognized that the propagation of the model (\ref{eq:eqm_psi}--\ref{eq:eqm_v_a}) fails due to the singularity, while the quaternion formulation 
(\ref{eq:q_eqm}--\ref{eq:q_va})
leads to correct numerical simulation results.
\section{Optimal control problem formulation}
\label{sec:optimal_control}

The optimization task is to control power cycles with optimal power output averaged
over complete cycles. 
This quantity is computed as integral over the mechanical power
\begin{equation}
  \bar{P} = \frac{1}{T} \int\limits_0^T \dot{l}\, F_{\rm tether} \;\mathrm{d}t =
  \frac{\varrho A C_{\rm R}}{2} \frac{E}{\sqrt{1\!+\!E^2}} \frac{1}{T} \int\limits_0^{T} \dot{l}\,
  v_{\rm a}^2\;\mathrm{d}t\label{eq:P_avg}
\end{equation}
It should be mentioned at this point that an 
{\em open loop} optimal control problem is set up, which is based on controls as model input, but does not involve any feedback control based on measured values.
\subsection{Setup of the optimization problem}
As steering speed will also be limited, the steering deflection $\delta$ is introduced as additional state to the model and the pod steering actuator speed, defined as $\dot{\delta}_{\rm s}$, is used as control, instead. 
Further, the state $W$ with $W(0)=0$ for computing the integral in (\ref{eq:P_avg}) is added.
Thus, the following relations have to be added to complete the set of EQM
(\ref{eq:q_eqm}--\ref{eq:q_va})
\begin{eqnarray}
  \dot{\delta} &=& \dot{\delta}_{\rm s} \label{eq:eqm_deltadot}\\
  \dot{W} &=& \dot{l} v_{\rm a}^2\label{eq:eqm_W}
\end{eqnarray}
The control vector now reads
\begin{equation}
  {\bf u} = \left[ \dot{\delta}_{\rm s}, v^{\rm (winch)} \right]^\top
\end{equation}
and the state 
\begin{equation}
  {\bf x} = \left[ W, \delta, l , q_0, q_1, q_2, q_3\right]^\top
\end{equation}
The NLP problem is formulated as minimization of
\begin{equation}
  f = -\frac{1}{T}\int\limits_0^{T}  \dot{l} v_{\rm
  a}^2 \;\mathrm{d}t + \epsilon_{\delta} \int\limits_0^{T} |\dot{\delta}(t) |^2 \;\mathrm{d}t
    = -\frac{W(T)}{T} + \epsilon_{\delta} \int\limits_0^{T}  |\dot{\delta}(t) |^2 \;\mathrm{d}t
\end{equation}
The second term introduces a penalty in order to achieve a smooth steering
behavior. The weighting factor $\epsilon_{\delta}$ has to be chosen appropriately.

The optimization, i.e.~minimization of $f$, has to be carried out 
by variation of $\vec{x}(t),\vec{u}(t),T$ subject to the
following constraints.
First the optimization is subject to the model equations of motion (\ref{eq:q_eqm}--\ref{eq:q_va}).
In addition, physical, geometrical and topological constraints have to be added, which will be discussed 
in the subsequent sections.
\subsection{Physical constraints}
Due to weight and geometrical design considerations for the airborne system,
the two limitations of limited steering speed and deflection of the control pod actuator have to be taken into account.
Hence, the steering speed limit is given by
  \begin{equation}
    | \dot{\delta}_{\rm s} | \leq \dot{\delta}_{\rm max} \label{eq:constraint_deltadot}
  \end{equation}
and the limited steering deflection is implemented by
  \begin{equation}
    | \delta | \leq \delta_{\rm max} \label{eq:constraint_delta}
  \end{equation}
respectively.
In addition, winching speed is limited
\begin{equation}
  v_{\rm min}^{\rm (winch)} \leq v^{\rm (winch)} \label{eq:constraint_v_winch}
\end{equation}
Note that this bound on the reel-in speed imposes implicitly an constraint on the windward position during the reel-in phase, as can be shown by computing the steady state angle \cite{Erhard2013b},\cite{Erhard2015a}. 
Further, in order to keep the system tethered, a
certain minimum air path speed and hence tether force is required.
Recalling (\ref{eq:eqm_v_a}) yields
\begin{equation}
  v_{\rm a} = E v_{\rm w} (q_0^2+q_1^2-q_2^2-q_3^2) - E \dot{l} \geq v_{\rm a,min} \label{eq:constraint_v_a}
\end{equation}
Note that this condition is needed to guarantee the validity of the model during the optimization process.
\subsection{Geometric constraints}
As already introduced, optimization is done for closed pumping orbits, which impose periodic boundary conditions to the optimization problem
\begin{eqnarray}
    \vec{q}(T) &=& \vec{q}(0) \label{eq:period_q}\\
    l(T) &=& l(0) \\
    \delta(T) &=& \delta(0)\label{eq:period_delta}
  \end{eqnarray}
We also add one non-periodic boundary condition in order to define the
initial value of the ``energy state'' $W$: 
\begin{equation}   
W(0) = 0 \label{eq:period_initial}
 \end{equation}
 The tether length has to be constrained in order to keep the cycle within a certain range of line length. This is accomplished by setting an upper bound to the state
\begin{equation}
  l \leq l_{\rm max}\label{eq:constraint_lmax}
\end{equation}
  
For a real world system, earth (or water) surfaces are a serious constraint on
the equations of motion. Practical safety considerations taking into
account flight trajectory deviations due to wind gusts recommend the
choice of a
certain minimum elevation angle $\theta_{\rm min}$.
Geometric considerations yield
\begin{equation}
  \frac{-\scalarproduct{\vec{r}}{\vec{e}_{z}}}
  {\scalarproduct{\vec{r}}{\vec{e}_{x}}}
  \geq \tan\theta_{\rm min}
\end{equation}
Using (\ref{eq:q2position}) results in the constraint
\begin{equation}
  (q_0^2+q_1^2-q_2^2-q_3^2)\tan\theta_{\rm min} + 2(q_1 q_3 - q_0 q_2)\leq 0\label{eq:constraint_thetamin}
\end{equation}  
It should be finally noted that the operational flight altitude strongly depends on the wind profile and ideally, the optimal flight altitude lies above the given safety limit. 
However, for the optimal control computations in this paper, the wind field is assumed constant and homogeneous for simplicity and therefore these constraints play a major role.
\subsection{Topological constraints}
\label{sec:topo_constraints}
A crucial point in the formulation of the optimization problem is the topology of the trajectory.
For operational reasons, basically in order to avoid twisting of the tether, the patterns of choice are lemniscates.
Hence, one wishes to preserve the topology of the pattern during the optimization process by imposing certain constraints. If this is not done, the optimizer tries to 'unwind' certain trajectory features, which in almost all cases leads to a failure of the optimization process.

The topological constraints aim at keeping the configuration of the
pattern, e.g.~two lemniscates as sketched in Fig.~\ref{fig:topo}, or
six, as in the actual computations of this paper.
\begin{vchfigure}
\includegraphics[width=12cm]{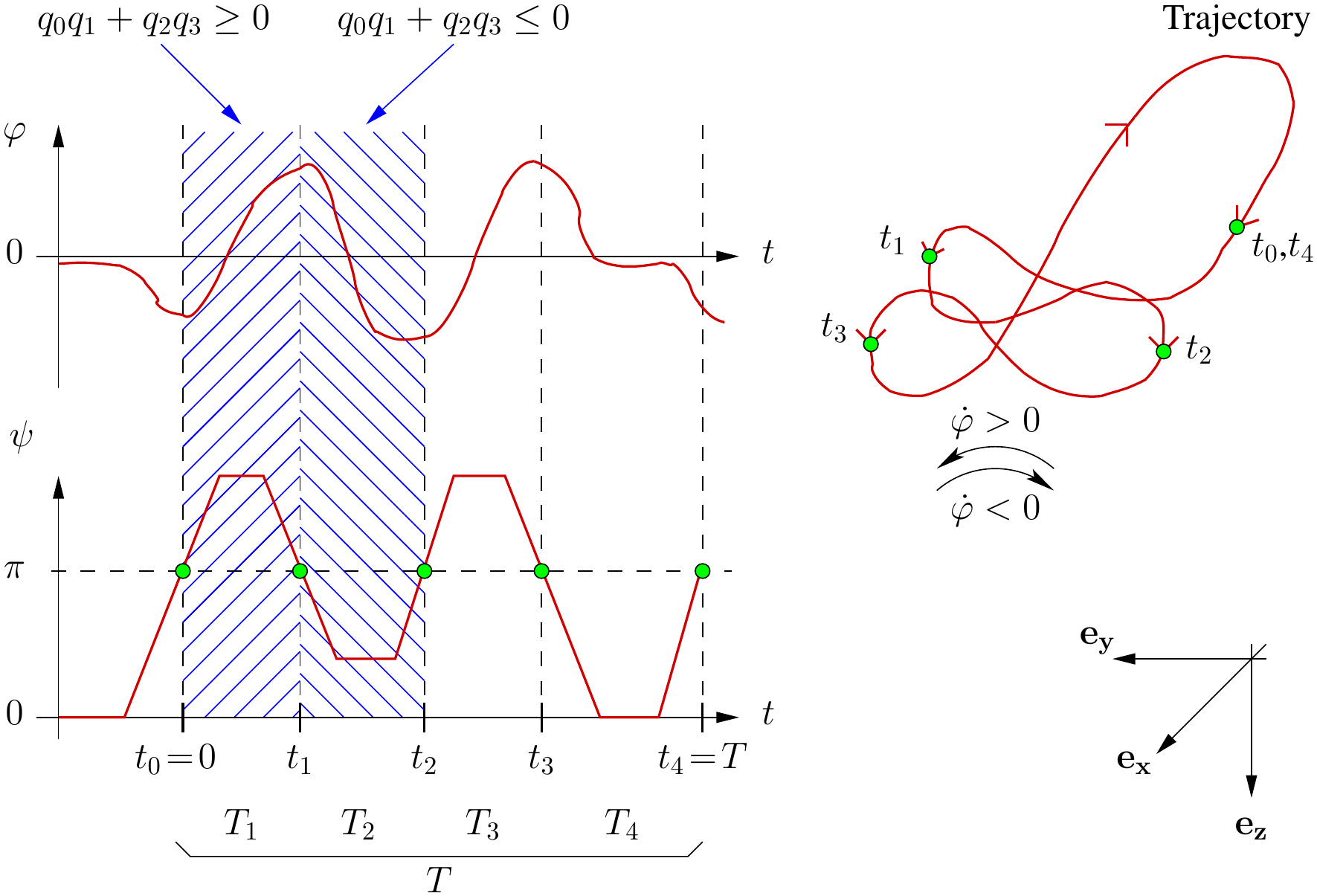}%
\vchcaption{Topological constraints shown for a pumping cycle consisting of two lemniscates. In order to avoid 'unwinding of the topology' by the optimizer, the trajectory horizon is divided into stages and additional constraints are imposed on those to preserve the topology.} 
\label{fig:topo}
\end{vchfigure}
The trajectory can be divided into stages, mainly parts featuring the property flying to the right (left), which corresponds to $\dot{\varphi}\!>\!0$ and $\dot{\varphi}\!<\!0$, respectively.
Formally, an even number of $N$ time intervals with duration $T_i$ is introduced.
Further, the time points $t_0, \dotsc ,t_N$ with the property $T_i=t_i - t_{i-1}$ are defined as shown in Fig.~\ref{fig:topo}. Note that $t_0=0$ and that $t_N=T$ is the overall pumping cycle time.

As $\dot{\varphi}$ is not a state of the system, (\ref{eq:eqm_phi})
can be used to formulate an equivalent constraint on $\psi$.
Considering the sign of $\sin\psi$ together with the behavior of the four-quadrant arctangent function $\mathrm{atan2}$ in (\ref{eq:q2varphi}), the following conditions can be deduced
\begin{equation}
  \begin{array}{c@{\hspace{1cm}}c@{\hspace{1cm}}c}
  q_0 q_3 - q_1 q_2 \geq 0 & t_{(2i)}\leq t \leq t_{(2i\!+\!1)} & i=1,\dotsc,(N/2) \\
  q_0 q_3 - q_1 q_2 \leq 0 & t_{(2i\!-\!1)}\leq t \leq t_{(2i)} & 
  \end{array}\label{eq:constraint_psi}
\end{equation}
The upper (lower) relation corresponds to 
$\dot{\varphi}\geq 0$ ($\dot{\varphi}\leq 0$).
%
%
%
\subsection{Summary of the optimal control problem}
The optimal control problem (OCP) can be stated formally as follows:
\begin{equation}\label{contocp-statement}
\!\!\!\!\!\!\!\!\!
\barr{r}\\
\stackrel{\mbox{minimize}}{\vec{x}(\cdot),\vec{u}(\cdot), t_1, \ldots, t_N}
\earr
  \quad E(\vec{x}(t_N),t_N)  + \int\limits_0^{t_N} L(\vec{u}(t)) \; \dd t
\end{equation}
\[
\left.\barr{rcll}
\mbox{subject to} \quad \quad \quad \quad \quad \\
\!\dot{\vec{x}}(t)\!-\!\Phi(\vec{x}(t),\vec{u}(t)) \!\!\!& = & 
0, \; \quad  \quad   t \in [t_0,t_N], \quad &\mbox{(ODE model)}\\
r_{\mathrm{boundary}}(\vec{x}(t_0), \vec{x}(t_N)) & =
&0,&\mbox{(boundary conditions)}\\
h(\vec{x}(t),\vec{u}(t))
&\leq& 0, \;   \quad \quad   t \in [t_{0},t_{N}],  &\mbox{(path constraints)}\\
h_i(\vec{x}(t))
&\leq& 0, \;   \quad \quad   t \in [t_{i-1},t_{i}], \quad i = 1, \ldots,
N.&\mbox{(stage wise path constraints)}
\earr
\right.
\]
with
\begin{equation}
  E(\vec{x}(t_N),t_N) = -\frac{W(t_N)}{t_N} 
  \quad\quad\mbox{and}\quad\quad
  L(\vec{u}(t))=\epsilon_{\delta} |\dot{\delta}(t)|^2
\end{equation}
The set of equations for the ODE model $\Phi(\cdot)$ is defined by (\ref{eq:eqm_quaternions}--\ref{eq:q_va}) and (\ref{eq:eqm_deltadot}--\ref{eq:eqm_W}).
The boundary conditions $r_{\rm boundary}(\cdot)$ are given by the
periodicity and initial state constraints (\ref{eq:period_q}--\ref{eq:period_initial}), the path constraints 
$h(\cdot)$ by (\ref{eq:constraint_deltadot}--\ref{eq:constraint_v_a}),
(\ref{eq:constraint_lmax}),
(\ref{eq:constraint_thetamin})
 and the stage wise path constraints $h_i(\cdot)$ by (\ref{eq:constraint_psi}).
\section{Numerical computation of optimal pumping orbits}
\label{sec:results}
In order to numerically solve the optimal control problem, 
a direct method is used to transfer
the continuous time OCP 
into a nonlinear programming problem (NLP).
This is accomplished by using the direct multiple shooting 
method. 
The direct multiple shooting method that was originally developed by
Bock and Plitt \cite{Bock1984} 
performs first a division of the time horizon into subintervals,
that we will call multiple shooting (MS) intervals in the sequel, and it uses 
piecewise control discretization on these intervals.
The continuous time dynamic system is 
transformed to a discrete time system by using an embedded ordinary
differential equation (ODE) solver to solve the ODE on each MS
interval individually.
Details of the setup of the discrete time control problem will be
summarized in the following subsections, 
and finally numerical results will be presented, that have been obtained by using the
optimization environment  CasADi~\cite{Andersson2012a,Andersson2013b} 
and the NLP solver IPOPT~\cite{Waechter2006}.
\subsection{Multiple shooting and control discretization grid}
The choice of time discretization grid size choice is a compromise between keeping the total number of optimization variables low and at the same time reproducing the dynamics of the controls and states in an adequate way. 
In order to represent controls and states of different timescales, nested grids are introduced as depicted in Fig.~\ref{fig:grid}.
\begin{vchfigure}
\includegraphics[width=10cm]{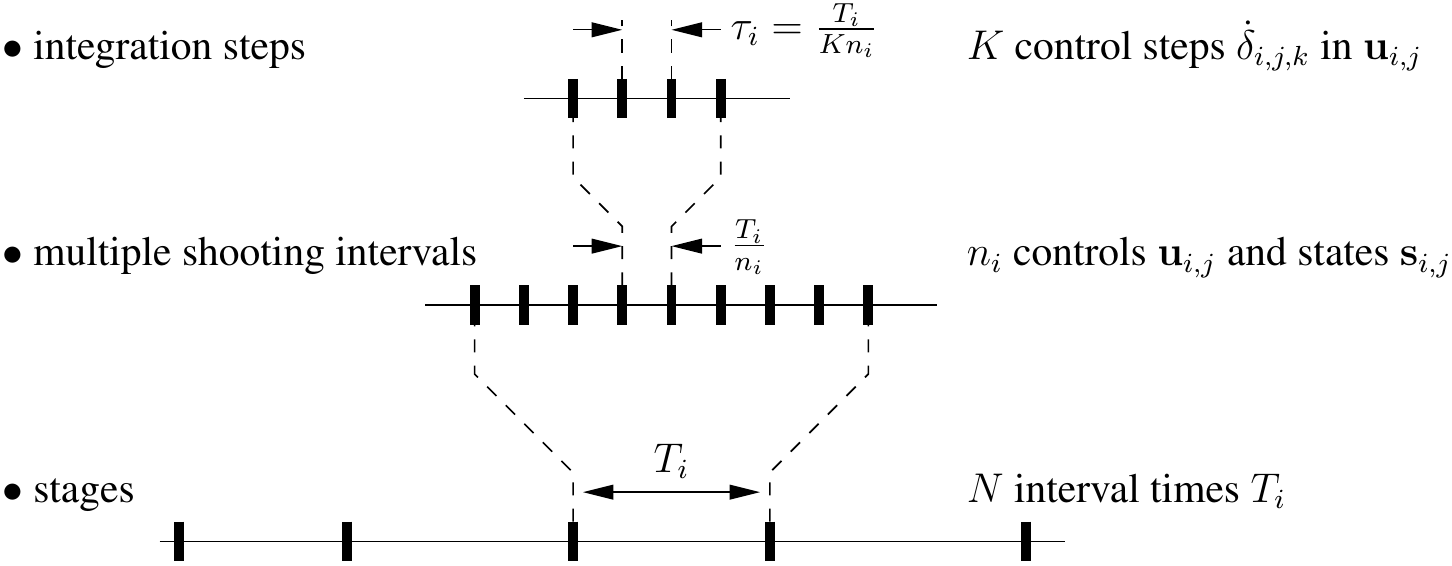}%
\vchcaption{Setup of the discretization grid.}
\label{fig:grid}
\end{vchfigure}
The upper level grid is the splitting of the trajectory into stages as introduced in Sect.~\ref{sec:topo_constraints}. The stage of duration $T_i$ is divided into 
$n_i$ equidistant multiple shooting (MS) intervals for the states and the winch control. 
Note that the stages can comprise different numbers $n_i$ of MS intervals in order to account for different phases as power and return phase, respectively.
These MS intervals are divided into $K$ subintervals for the kite steering actuator speed (control) to allow for a proper consideration
of the constraints (\ref{eq:constraint_deltadot}) in the solution.
The initial value on teach the control vector on each MS interval consists of $K+1$
components, i.e.
\begin{equation}
\vec{u}_{i,j} = \left[
  \dot{\delta}_{i,j,1},\dotsc ,\dot{\delta}_{i,j,K},v^{\rm (winch)}_{i,j}
  \right]^\top.
\end{equation}
The idea of multiple shooting is to solve the ODE on small intervals 
with duration $\frac{T_i}{n_i}$ starting with initial values
$\vec{s}_{i,j}$ 
 The used scheme is sketched in Fig.~\ref{fig:shooting}. 
\begin{vchfigure}
\includegraphics[width=10cm]{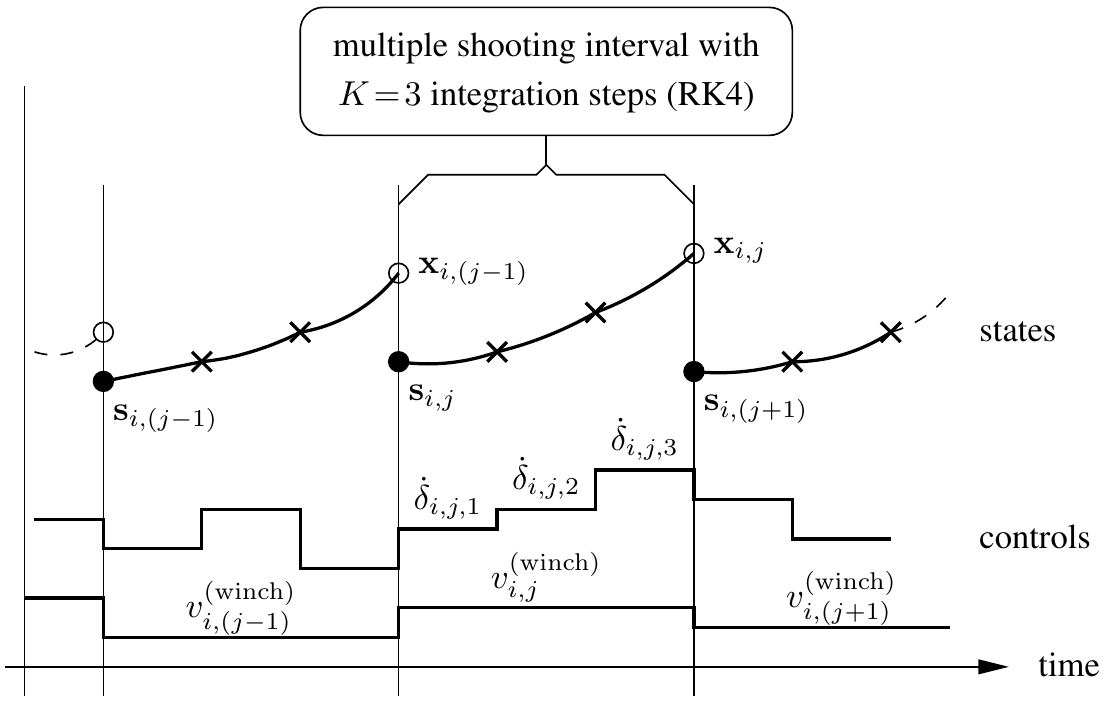}%
\vchcaption{Principle of the implemented shooting scheme for $K\!=\!3$.
The ODE model is solved on intervals by $K$ integration 
steps, starting at $\vec{s}_{i,j}$ and ending at $\vec{x}_{i,j}$. 
The trajectory is 'connected' by imposing constraints, e.g.~$\vec{s}_{i,(j+1)}-\vec{x}_{i,j}=0$. 
Note that in contrast to the usual implementation of multiple shooting, 
several piecewise constant control steps are taken per MS interval for
one of the controls.
}
\label{fig:shooting}
\end{vchfigure}
The integration is performed 
from an initial state $\vec{s}_{i,j}$ 
to the final state $\vec{x}_{i,j}$, using the controls
$\vec{u}_{i,j}$. Here, we are  utilizing
$K$ integration steps of the classical Runge-Kutta integrator of order
4 (RK4), each with a step size $\tau_{i} = \frac{T_i}{K n_i}$. 
In order to give a formal representation,
the single RK4 integration step  
$\Phi_{\rm RK}: \vec{s} \rightarrow \vec{x}$
shall be defined as 
\begin{equation}
  \vec{x} = \Phi_{\rm RK}(\tau_i, v^{\rm (winch)},\dot{\delta},\vec{s})
\end{equation}
The result of the $K$ RK4 integration steps on the interval with duration $\frac{T_i}{n_i}$ is then given by
\begin{equation}
  \vec{x}_{i,j}(s_{i,j}, u_{i,j}) \doteq
    \underbrace{\Phi_{\rm RK}(\tau_{i},v^{\rm (winch)}_{i,j},\dot{\delta}_{i,j,K},
      \underbrace{\Phi_{\rm RK}(\tau_{i},v^{\rm (winch)}_{i,j},\dot{\delta}_{i,j,(K-1)},\dotsc,
        \underbrace{\Phi_{\rm RK}(\tau_{i},v^{\rm (winch)}_{i,j},\dot{\delta}_{i,j,1},\vec{s}_{i,j})}_{\vdots} \cdots
      )}
    )}_{\mbox{$K$ integration steps}} 
\end{equation}
It should be emphasized that separate piecewise constant control 
steps $\dot{\delta}_{i,j,k}$ on the subgrid are provided for the integration steps of the MS interval.

In summary, the
optimization vector $\vec{w}=[\vec{w}_1^\top, \ldots,
\vec{w}_N^\top]^\top$ 
can be subdivided into subvectors with all variables $\vec{w}_i$
belonging to stage $i$, which are 
given by 
\begin{equation}
\vec{w}_i=  \left[\vec{s}_{i,0}^\top, \vec{u}_{i,0}^\top, \ldots,
  \vec{s}_{i,(n_1-1)}^\top, \vec{u}_{i,(n_1-1)}^\top,
  \vec{s}_{i,n_i}^\top, T_i\right]^\top.
\end{equation}

\subsection{Summary of the discretized optimal control problem}
The complete discretized OCP can be stated as follows
\begin{equation}
\begin{array}{cc}
  \begin{array}{c}\\{\mbox{minimize}}\\{\vec{w}}\end{array}
  & E(\vec{s}_{N, n_N}, \sum_{i=1}^N T_i) + \sum_{i=1}^N
  \sum_{j=0}^{n_i-1}  \ell_i(\vec{u}_{i,j}, T_i)
  \end{array}
\end{equation}
    
\begin{equation}
  \begin{array}{r@{\hspace{.7cm}}l@{\hspace{.7cm}}l}
  \multicolumn{2}{l}{\mbox{subject to}} \\
  \vec{x}_{i,j}(\vec{s}_{i,j}, \vec{u}_{i,j}) - \vec{s}_{i,(j+1)} = 0 & i
  = 1,\dotsc ,N; j = 0, \dotsc ,(n_i\!-\!1) &  \mbox{(continuity
    between MS intervals)}\\
  \vec{s}_{i ,n_i} - \vec{s}_{(i+1) ,0} = 0 & i = 1,\dotsc ,(N\!-\!1)
  &  \mbox{(continuity between stages)}\\
  r_{\rm boundary} (\vec{s}_{1,0}, \vec{s}_{N, n_N}) = 0 & 
  & \mbox{(boundary conditions)} \\
  \tilde{h}(\vec{s}_{i,j},\vec{u}_{i,j}) \leq 0 & i
  = 1,\dotsc ,N; j = 0, \dotsc ,(n_i\!-\!1) & \mbox{(path constraints)} \\
  \tilde{h}_i(\vec{s}_{i,j}) \leq 0 & i = 1,\dotsc ,N ; j = 0, \dotsc ,(n_i\!-\!1)  & \mbox{(stagewise path constraints)}
  \end{array} \nonumber
\end{equation}
with
\begin{equation}
  \ell_i(\vec{u}_{i,j},T_i) = 
  \left(
  \frac{T_i}{n_i} \sum\limits_{k=1}^K\epsilon_{\delta} \dot{\delta}_{i,j,k}^2 
  \right) +
  \frac{T_i}{n_i} \epsilon_v\left(v^{\rm (winch)}_{i,j}-v^{\rm (winch)}_{{\rm next}(i,j)}\right)^2
\end{equation}
where next$(i,j)$ denotes the subsequent index pair in time. Note that
the second term penalizes changes in winch speed (i.e.~acceleration)
and has been added to the discretized OCP only.  The discretized path constraints 
$\tilde{h}(\cdot)$ evaluate the constraints (\ref{eq:constraint_deltadot}--\ref{eq:constraint_v_a}),
(\ref{eq:constraint_lmax}),
(\ref{eq:constraint_thetamin}) on the multiple shooting grid, 
 and the discretized stage wise path constraints $\tilde{h}_i(\cdot)$
 evaluate (\ref{eq:constraint_psi}) on the grid.

\subsection{Numerical results}
\begin{figure}
\includegraphics[width=\textwidth]{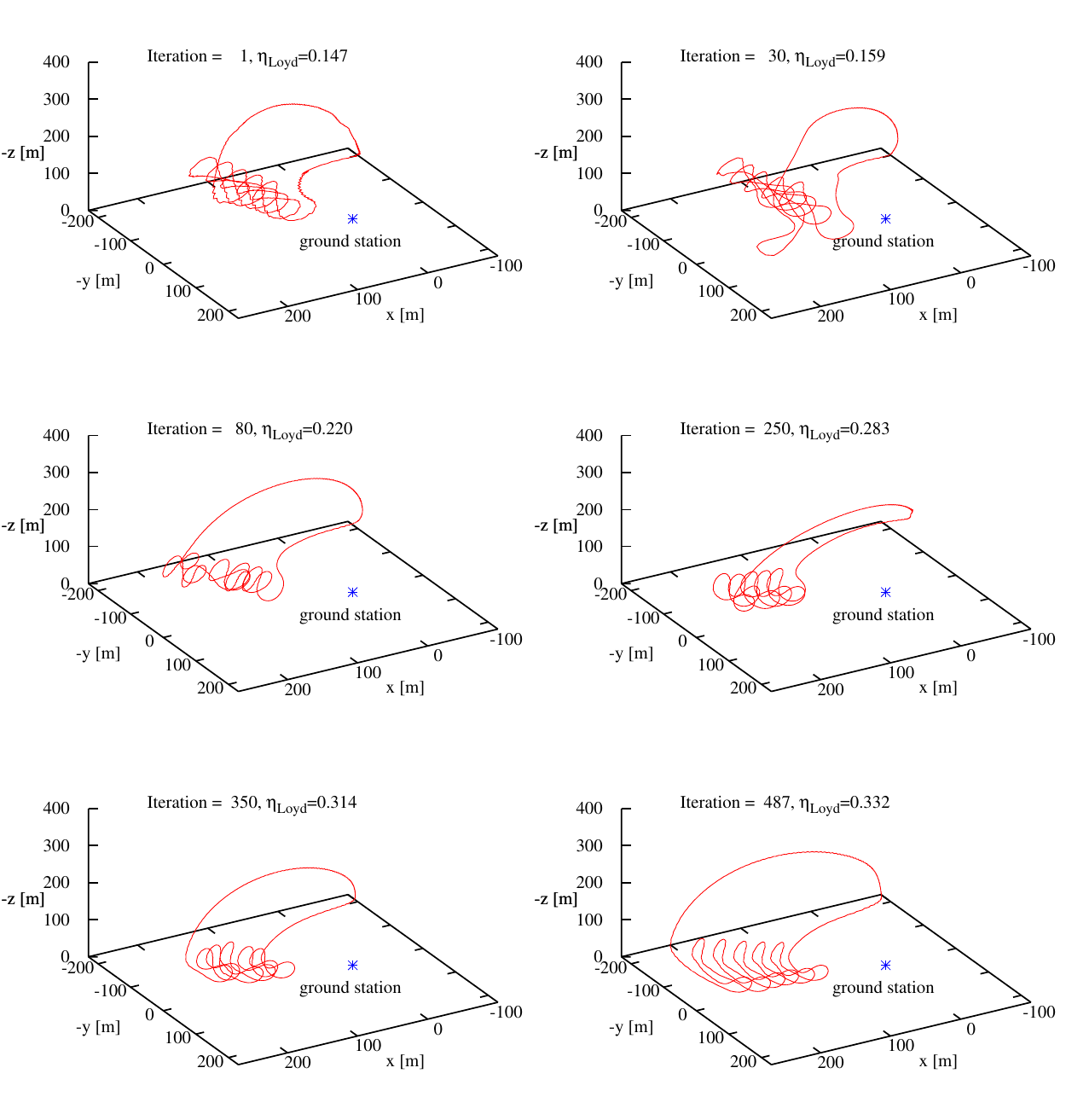}%
\caption{Evolution of pumping cycles during the optimization process. Starting with the experimentally flown trajectory (top left), the 3d trajectories for different iteration steps (number given in figure) to the optimal solution (bottom right) are shown.}
\label{fig:result_trajectories}
\end{figure}

In the following, numerical results for optimal complete pumping cycle will be presented. 
The parameters are given in Table \ref{tab:parameters} and are related to the SkySails functional prototype of Fig.~\ref{fig:sks_setup}.
As initial guess, an experimentally flown trajectory is used. The topology of 6 lemniscates leads to $N\!=\!12$ and the (initial) cycle time of $T\!\approx\!170$\,s is divided up into $\sum_{i=1}^N n_i=250$ intervals. The number of RK4 integration steps is chosen as $K\!=\!3$ on each MS integration interval.
Thus, the optimization vector contains a total of 2762 optimization variables.
The 3d trajectories at different iteration steps are shown in 
Fig.~\ref{fig:result_trajectories}.
It can observed that optimization starts far from optimum and runs through some 
crude path configurations towards the optimal solution.
This behavior could be attributed to the interior point method of the used IPOPT solver.

The optimization figures of merit are given by the Loyd factor $\eta_{\rm
Loyd}$, which is defined by 
\begin{equation}
  \eta_{\rm Loyd} = \frac{\bar{P}}{P_{\rm Loyd}}
\end{equation}
This factor corresponds to the Loyd limit, which corresponds to the generated power $P_{\rm Loyd}$ for a 
kite flying 'maximal crosswind' (at $\vartheta\!=\!0$) continuously, i.e.~no retraction
phase, compare (\ref{loydsformula}) and \cite{Loyd1980}. 
For the derived model, this quantity is given by
\begin{equation}
  P_{\rm Loyd} = \frac{\varrho C_{\rm R}A}{2} \frac{4E^2}{27} 
  \frac{E}{\sqrt{1\!+\! E^2}} v_{\rm w}^3
\end{equation}

\begin{figure}
\includegraphics[width=\textwidth]{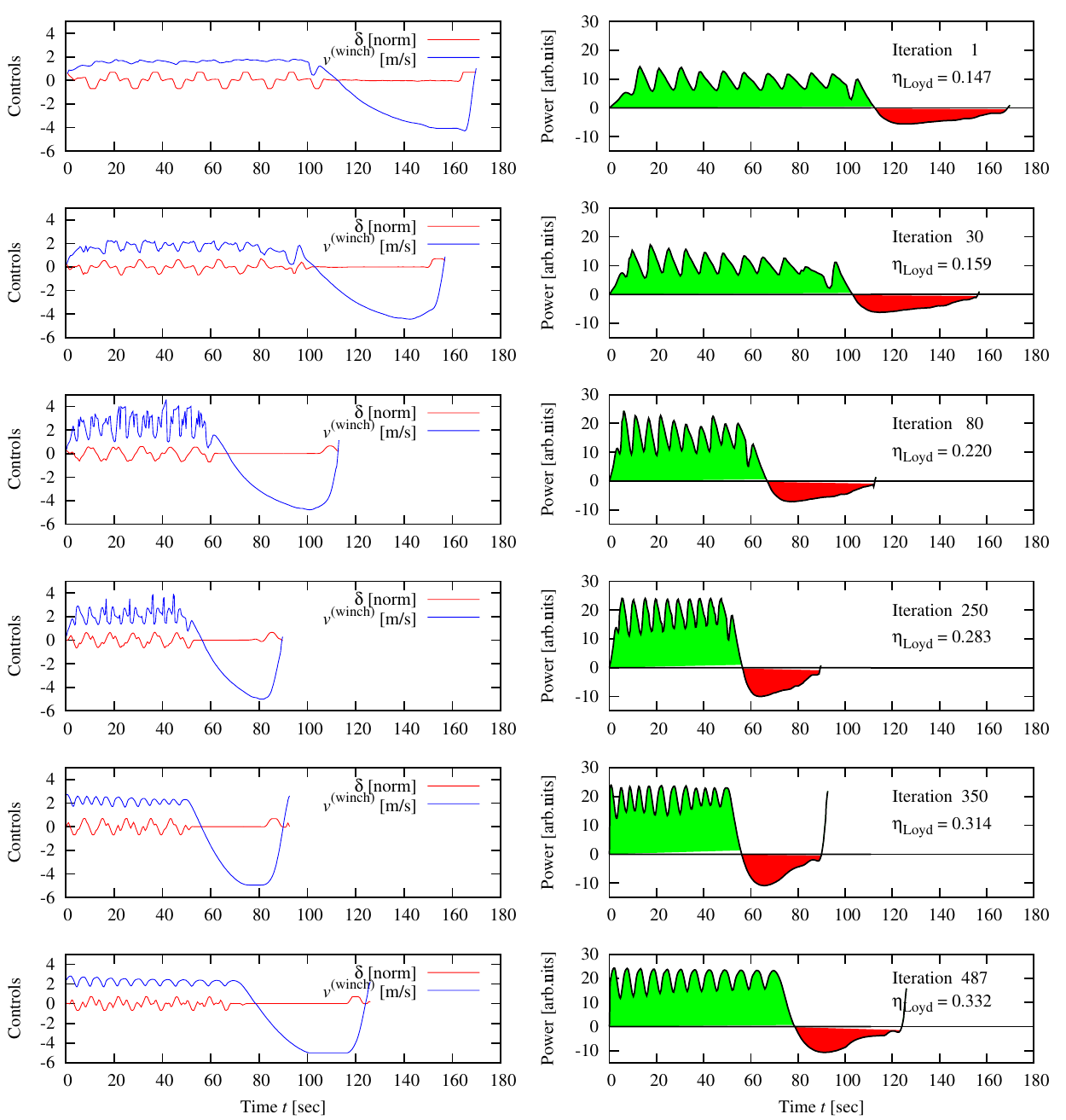}%
\caption{Time series for controls $\delta$ and $v^{\rm (winch)}$ and power, computed as tether force times reeling speed $F_{\rm tether} \dot{l}$. The upper row shows the experimentally flown pumping cycle, which is used initial guess. The lower row gives the optimal solution. Iteration numbers for the intermediate steps are annotated in the respective subplots. }
\label{fig:result_timeseries}
\end{figure}

The according time series for the controls and power, given as $F_{\rm tether} \dot{l}$, are presented in Fig.~\ref{fig:result_timeseries}.

In summary, it can be stated that the optimizer increased the power generation of the pumping cycle by more than a factor of two from $\eta_{\rm Loyd} \!\approx\! 0.15$ to $\eta_{\rm Loyd} \!\approx\! 0.33$.
This is partly due to the fact that the optimized figure lies deeper in the crosswind region. Further, transfer and return phases are carried out more efficiently.
The lemniscates are deformed to minimal curves fulfilling the imposed constraints.
For sake of a fair assessment, it should be noted that the experimental trajectory 
was taken from a test flight aiming at flight controller development rather than  optimal power output.
Nevertheless, valuable suggestions for improvement of even these 'classical' control setups can be extracted from such optimization results, see e.g.~the winch controller implementation for the transfer phase in \cite{Erhard2015a}.
\section{Conclusion}
\label{sec:conclusion} 
In conclusion, a singularity free model for tethered kite dynamics based on quaternions has been derived, which can be broadly applied for simulation and optimization purposes.
It has been demonstrated that optimization runs of complete pumping
cycles within a few minutes are feasible with this model.
This opens up the way to further extended and systematic studies of AWE efficiencies.
It should be finally remarked that these might most likely involve further extensions to the model to take into account e.g.~mass effects, tether drag and a wind profile.
\begin{acknowledgement}
  We gratefully thank Mario Zanon for valuable discussions on the
  implementation of the OCP. This research was supported by KUL PFV/10/002 OPTEC;
  Eurostars SMART; IUAP P7 (DYSCO); EU:  FP7-TEMPO (MCITN-607957), 
H2020-ITN AWESCO (642682) and by the
ERC Starting Grant HIGHWIND (259166).
\end{acknowledgement}

%
%
 
%
\providecommand{\WileyBibTextsc}{}
\let\textsc\WileyBibTextsc
\providecommand{\othercit}{}
\providecommand{\jr}[1]{#1}
\providecommand{\etal}{~et~al.}

\end{document}